\documentclass[a4paper,12pt]{article}

\usepackage[latin1]{inputenc}

\usepackage[T1]{fontenc}

\usepackage{graphicx}                      

\usepackage[all]{xy}
\usepackage{amsmath,amsfonts,amssymb}

\usepackage{mathptmx}

\input{macro}
\newcommand{\ff}{\emph{focus-focus}}
\renewcommand{\F}{\mathcal{F}}

\title{On semi-global invariants for \ff\ singularities} 

\author{V\~u Ng\d oc San}

\date{june 2001}

\begin{document}
\maketitle

\vskip 1cm
\begin{center} 
  Institut Fourier\\
  Unit{\'e} mixte de recherche CNRS-UJF 5582\\
  BP 74, 38402-Saint Martin d'H\`eres Cedex (France).\\
  \emph{Email~:~\texttt{San.Vu-Ngoc@ujf-grenoble.fr}}
\end{center}
\vskip 1cm

\begin{abstract}
  This article gives a classification, up to symplectic equivalence,
  of singular Lagrangian foliations given by a completely integrable
  system of a 4-dimensional symplectic manifold, in a full
  neighbourhood of a singular leaf of \ff\ type.
\end{abstract}

\vfill

\noindent\textbf{Keywords~: } symplectic geometry, Lagrangian fibration,
completely integrable systems, Morse singularity, normal forms.

\noindent\textbf{Math. Class.~:}
37J35, 
37J15, 
70H06, 
37G20, 
53D22, 
53D12, 
70H15. 

\newpage

\section{Introduction}

In the study of completely integrable Hamiltonian systems, and more
generally for any dynamical system, finding normal forms is often the
easiest way of understanding the behaviour of the trajectories. Normal
forms generally deal with a \emph{local} issue. But the locality here
depends on one's viewpoint: one can be local near a point, an orbit,
or any invariant submanifold. If $F=(H_1,\dots,H_n)$ is a completely
integrable system on a $2n$-symplectic manifold $M$ (meaning that
$\{H_j,H_i\}=0$), several normal forms hold:
\begin{itemize}
\item near a point $m$ where $dH_j(m)$, $j=1,\dots,n$ are linearly
  independent, one can construct \emph{Darboux-Carathéodory}
  coordinates: a neighbourhood of $m$ is symplectomorphic to a
  neighbourhood of the origin in $\RM^{2n}$ with its canonical
  coordinates $(x,\xi)$, in such a way that $H_j-H_j(m)=\xi_j$.
\item if $c$ is a regular value of $F$, one has near any compact
  connected component $\Lambda_c$ of $F^{-1}(c)$ the
  \emph{Liouville-Arnold} theorem which states that the system is
  symplectomorphic to a neighbourhood of the zero section of
  $T^*(\T^n)$ in such a way that there is a change of coordinates
  $\Phi$ in $\RM^n$ such that $F\circ\Phi=(\xi_1,\dots,\xi_n)$. Here
  $\T^n$ is the torus $\RM^n/2\pi\ZM^n$ and the cotangent bundle
  $T^*(\T^n)$ is equipped with canonical coordinates $(x,\xi)$.
\end{itemize}
The first one is typically a local normal form, while I would refer to
the Liouville-Arnold theorem as a \emph{semi-global} result, for it
classifies a neighbourhood of a whole invariant Lagrangian leaf
$\Lambda_c$. These two statements above are now fairly standard. They
can be extended in different directions: a) trying to globalise: what
can be said at the level of the whole fibration of regular fibres
$\Lambda_c$ ? This of course involves more topological invariants, as
described in Duistermaat's paper \cite{duistermaat}; b) including
critical points, which is the main incentive for this article.

A Morse-Bott like theoretical study of critical point of completely
integrable Hamiltonian systems exists, which yields a local symplectic
classification of non-degenerate singularities (see Eliasson
\cite{eliasson-these}). These results have been used by Nguyên Tiên
Zung \cite{zung} (extending previous results by Fomenko) to obtain a
topological semi-global classification of the singular foliation. This
work does not give the corresponding smooth symplectic classification,
where new semi-global invariants show up, as demonstrated in the
``1-D'' (one degree of freedom, i.e. $n=1$) case by
\cite{dufour-mol-toul}. The point of our present article is to extend
the results of \cite{dufour-mol-toul} to the 2-D case of \ff\ 
singularities. Note that our arguments could easily be applied in the
1-D case, thus supplying for the lack of proofs in
\cite{dufour-mol-toul}.

Between the pure topological classification of the singular foliation
and the ``exact'' symplectic classification, some other interesting
notions of equivalence have been introduced (see eg.
\cite{bolsinov-fomenko}), which are all weaker than what we shall
present here.

The semi-global viewpoint seems to be able to shed some new light in
semiclassical mechanics, where a quantum state is associated to a
Lagrangian submanifold. Quantum states associated to singular
manifolds have a particularly rich structure, strongly linked to the
local (for this, see \cite{san-fn}) and semi-global symplectic
invariants of the foliation. We expect to return on this in a future
paper.

\section{Statement of the result}

In this article, $(M,\omega)$ is a 4-dimensional symplectic manifold,
equipped with the symplectic Poisson bracket $\{\cdot,\cdot\}$. Any
smooth function $H$ on $M$ gives rise to a Hamiltonian vector field
denoted by $\ham{H}$.

The word smooth always means of $\Cinf$ category and a function $f$ is
said \emph{flat} at a point $m$ if $f$ and all its derivatives vanish
at $m$.

\begin{defi}
  A map $F=(H_1,H_2)$ defined on some open subset $U$ of $M$ with
  values in $\RM^2$ is called a \textbf{momentum map} if $dF$ is
  surjective almost everywhere in $U$ and $\{H_1,H_2\}=0$.
\end{defi}
\begin{defi}
  A \textbf{singular Liouville foliation} $\F$ is a disjoint union of
  connected subsets of $M$ called leaves for which there exists a
  momentum map $F$ defined in some neighbourhood $\Omega$ of $\F$ such
  that the leaves of $\F$ are the connected components of the level
  sets $F^{-1}(c)$, for $c$ in some open subset of $\RM^2$.
\end{defi}

The total space of the foliation is also denoted by $\F$. The above
definition implies that $\F$ is an open subset of $M$.

\begin{defi}
  Let $m\in\F$. The maximum of the set $\{\textup{rank}(dF(m)), F
  \textrm{ defining } \F \}$ is called the rank of $m$.  $m$ is called
  \textbf{regular} if its rank is maximal ($=2$). Otherwise it is
  called \textbf{singular}.
\end{defi}
 
If $m$ is a regular point, then there is an open neighbourhood of $m$
in which all points are regular, and if $F_1$ $F_2$ are associated
momentum maps near $m$, one has $F_1=\phy\circ F_2$, for some local
diffeomorphism $\phy$ of $\RM^2$ (these facts come from the local
submersion theorem). 

Note that the condition $\{H_1,H_2\}=0$ implies that the leaves are
local Lagrangian manifolds near any regular point. However, the
foliation near a regular leaf (=a leaf without any singular point) is
not the most general Lagrangian foliation (which would be defined as a
foliation admitting \emph{locally} associated momentum maps), since
the latter does not necessarily admit a global momentum map (see
\cite{weinstein-symplectic}).

In what follows, the word ``Liouville'' is often omitted.  If
$m\in\F$, we denote by $\F_m$ the leaf containing $m$.

\begin{defi}
  \label{defi:focus}
  A singular Liouville foliation $\F$ is called of \textbf{simple \ff\ 
    type} whenever the following conditions are satisfied:
  \begin{enumerate}
  \item $\F$ has a unique singular point $m$;
  \item the singularity at $m$ is of \ff\ type;
  \item the leaf $\F_m$ is compact.
  \end{enumerate}
  The leaf $\F_m$ is called the \ff\ leaf.
\end{defi}

Recall that the second condition means that there exists a momentum
map $F=(H_1,H_2)$ for the foliation at $m$ such that the Hessians of
$H_1$ and $H_2$ span a subalgebra of quadratic forms that admits, in
some symplectic coordinates $(x,y,\xi,\eta)$, the following basis:
\begin{equation}
q_1=x\xi+y\eta,\qquad q_2=x\eta-y\xi.\label{equ:cartan}
\end{equation}
This implies that \ff\ points are isolated, which ensures that the
above definition is non-void. Note that \ff\ singularities are one of
the four types of singularities of Morse-Bott type in dimension 4, in
the sense of Eliasson \cite{eliasson}.

\begin{defi}
  Two singular foliations $\F$ and $\tilde{\F}$ in the symplectic
  manifolds $(M,\omega)$ and $(\tilde{M},\tilde{\omega})$ are
  \textbf{equivalent} is there exists a symplectomorphism
  $\phy:\F\fleche\tilde{\F}$ that sends leaves to leaves.
\end{defi}

\begin{defi}
  Let $\F$ and $\tilde{\F}$ be singular foliations in $M$, and
  $m\in\F\cap\tilde{\F}$ such that $\F_m=\tilde{\F}_m$. The
  \textbf{germs} of $\F$ and $\tilde{\F}$ at $\F_m$ are equal if and
  only if there exists a saturated neighbourhood $\Omega$ of $\F_m$ in
  $\F$ such that $\F\cap\Omega=\tilde{\F}\cap\Omega$.
\end{defi}

The classification of germs of Liouville foliations near a compact
regular leaf is given by the Liouville-Arnold theorem that asserts
that they are all equivalent to the horizontal fibration by tori of
$T^*\T^n$. The presence of singularities imposes more rigidity, and we
have the following theorem (which is natural in view of
\cite{dufour-mol-toul}):

\begin{theo}
  \label{theo:abstrait}
  The set of equivalence classes of germs of singular Liouville
  foliations of \ff\ type at the \ff\ leaf is in natural bijection
  with $\RM\formel{X,Y}_0$, where $\RM\formel{X,Y}$ is the algebra of
  real formal power series in two variables, and $\RM\formel{X,Y}_0$
  is the subspace of such series with vanishing constant term.
\end{theo}

This formal statement does not contain the most interesting part of
the result, which is the geometric description of the power series
involved (it is essentially the Taylor series of a regularisation of
some action integral). The rest of the paper is devoted to this
description -- which is the ``$\Rightarrow$'' sense of the theorem,
and to the proof of the ``$\Leftarrow$'' sense, for which we provide a
normal form corresponding to any given power series in
$\RM\formel{X,Y}_0$.

The articles ends up with a sketchy argument as to how the result can
be extended to handle the case of several \ff\ points in the singular
leaf.

\section{The regularised action}

\begin{figure}[htbp]
  \begin{center}
    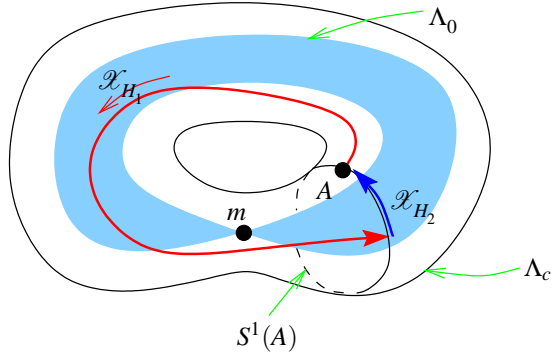
    \caption{Construction of the ``periods'' $\tau_j(c)$}
    \label{fig:monodromy}
  \end{center}
\end{figure}

Let $\F$ be a singular foliation of simple \ff\ type. Then in some
neighbourhood $U$ of the \ff\ point $m$, the following linearisation
result holds (Eliasson \cite{eliasson-these}): there exist symplectic
coordinates in $U$ in which the map $(q_1,q_2)$ (defined
in~(\ref{equ:cartan})) is a momentum map for the foliation.  Notice
therefore that, contrary to what the picture of
Figure~\ref{fig:monodromy} may suggest, $\F_m$ is diffeomorphic near
$m$ to the union of two 2-dimensional planes transversally
intersecting at $m$. Let $A$ be a point in $\F_m\cap U\setminus\{m\}$,
and $\Sigma$ be a small 2-dimensional surface transversal to the
foliation at $A$, and $\Omega$ be the open neighbourhood of $\F_m$
consisting of leaves intersecting $\Sigma$. In what follows, we
restrict the foliation to $\Omega$.

Let $\tilde{F}$ be a momentum map for the whole foliation $\F$
satisfying the hypothesis of Definition~\ref{defi:focus}. In a
neighbourhood of $\Sigma$, $\tilde{F}$ and $q=(q_1,q_2)$ are regular
local momentum maps, hence $q=\phy\circ\tilde{F}$, for some local
diffeomorphism $\phy$ of $\RM^2$. Now let $F=\phy\circ\tilde{F}$. It
is a global momentum map for $\F$ that extends $q$. We denote
$F=(H_1,H_2)$ and $\Lambda_c=F^{-1}(c)$.

Near $m$, the Hamiltonian flow of $q_2$ is $2\pi$-periodic, and --
assuming $U$ to be invariant with respect to this flow -- the
associated $S^1$-action is free in $U\setminus\{m\}$. Since this
action commutes with the flow of $H_1$, the $H_2$-orbits must be
periodic of primitive period $2\pi$ for any point in a (non-trivial)
trajectory of $\ham{H_1}$. On the leaf $\F_m=\Lambda_0$, these
trajectories are homoclinic orbits for the point $m$, which implies
that the flow of $H_2$ generates an $S^1$-action on a whole
neighbourhood of $\F_m$ (see \cite{san-focus} for details).

For any point $A\in\Lambda_c$, $c$ a regular value of $F$, let
$\tau_1(c)>0$ be the time of first return for the $\ham{H_1}$-flow to
the $\ham{H_2}$-orbit through $A$, and $\tau_2(c)\in\RM/2\pi\ZM$ the
time it takes to close up this trajectory under the flow of
$\ham{H_2}$ (see Fig.~\ref{fig:monodromy}). These times are
independent of the initial point $A$ on $\Lambda_c$.

For any regular value $c$ of $F$, the set of points $(a,b)\in\RM^2$
such that $a\ham{H_1}+b\ham{H_2}$ has a $1$-periodic flow on
$\Lambda_c$ is a sublattice of $\RM^2$ called the \emph{period
  lattice} \cite{duistermaat}.  The vector fields
$\tau_1\ham{H_1}+\tau_2\ham{H_2}$ and $2\pi\ham{H_2}$ both define
$1$-periodic flows, hence $(\tau_1,\tau_2)$ and $(0,2\pi)$ form a
$\ZM$-basis of the period lattice (see Remark~\ref{rema:naturality}).
As we shall see, the classification we are looking for relies on the
behaviour of this basis as $c$ tends to $0$. One immediate fact is
that the cycle associated to $\ham{H_2}$ shrinks to a point (vanishing
cycle). On the other hand, the coefficients of the first vector field
display a logarithmic divergence, as stated in the following
proposition.

\begin{prop}
  \label{prop:period}
  Let $\ln c$ be the some determination of the complex logarithm,
  where $c=(c_1,c_2)$ is identified with $c_1+ic_2$. Then the
  following quantities
  \[
  \left\{
    \begin{array}{ccl}
      \sigma_1(c) & = & \tau_1(c)+\Re(\ln c) \\
      \sigma_2(c) & = & \tau_2(c)-\Im(\ln c)
    \end{array}
  \right.
  \]
  extend to smooth and single-valued functions in a neighbourhood of
  $0$.  The differential 1-form
  \[
  \sigma:=\sigma_1dc_1+\sigma_2dc_2
  \]
  is closed.
\end{prop}
\begin{demo}
  As before, let $U$ be the neighbourhood of $m$ found using
  Eliasson's result, with canonical coordinates $(x,y,\xi,\eta)$.  In
  $U$, we use the complex coordinates $z=(z_1,z_2)$ with $z_1=x+iy$
  and $z_2=\xi+i\eta$, so that $q_1(z)+iq_2(z)=\bar{z_1}z_2$. The flow
  of $q_1$ is
  \begin{equation}
    (z_1(t),z_2(t))=(e^t z_1(0)),e^{-t} z_2(0),
    \label{equ:flowq1}
  \end{equation}
  while the flow of $q_2$ is the $S^1$-action given by
  \begin{equation}
    (z_1(t),z_2(t))=e^{it}(z_1(0),z_2(0)).
    \label{equ:flowq2}
  \end{equation}
  
  Fix some small $\epsilon>0$. Then the local submanifolds
  $\Sigma_u=\{z_1=\epsilon, |z_2| \textrm{ small}\}$ and
  $\Sigma_s=\{|z_1| \textrm{ small}, z_2=\epsilon,\}$ are transversal
  to the foliation $\Lambda_c=\{(z_1,z_2),\quad \bar{z_1}z_2=c\}$;
  therefore, the intersections $A(c):=\Sigma_u\cap\Lambda_c$ and
  $B(c):=\Sigma_s\cap\Lambda_c$ are smooth families of points.
  
  The $S^1$-orbits of $\Sigma_{u/s}$ form two small hypersurfaces
  transversal to the flow of $q_1$; therefore one can uniquely define
  $\tau_1^{A,B}(c)$ as the time of first hit on $\Sigma_s$ for the
  $\ham{H_1}$-flow starting at $A(c)$ (and hence flowing outside of
  $U$), and $\tau_2^{A,B}(c)$ as the time it takes to finally reach
  $B(c)$ under the $\ham{H_2}$-flow.  $\tau_1^{A,B}(c)$ and
  $\tau_2^{A,B}(c)$ are smooth functions of $c$ in a neighbourhood of
  $0$.
  
  Interchanging the roles of $A$ and $B$ -- and thus of $\Sigma_u$ and
  $\Sigma_s$, the times $\tau_j^{B,A}(c)$ for $j=1,2$ are defined in
  the same way.  But since the corresponding flows now take place
  inside $U$, where a singular point occur, $\tau_j^{B,A}(c)$ is not
  defined for $c=0$. On the other hand, equations~\eqref{equ:flowq1}
  and \eqref{equ:flowq2} yield the following explicit formula:
  \begin{equation}
    \label{equ:tau}
    \tau_1^{B,A}(c)+i\tau_2^{B,A}(c) = \ln \frac{z_1(A)}{z_1(B)} = \ln
    z_1(A)\bar{z_2}(B) - \ln \bar{c} = \ln \epsilon^2 - \ln \bar{c}
  \end{equation}
  
  Writing now
  \[
  \tau_1(c)+i\tau_2(c) = \left(\tau_1^{A,B}(c)+\tau_1^{B,A}(c)\right)
  + i\left(\tau_2^{A,B}(c)+\tau_2^{B,A}(c)\right),
  \]
  using (\ref{equ:tau}), and the fact that $\ln \bar{c}=\ln|c|-i\arg
  c$, we obtain that
  \[
  \sigma_1(c)+i\sigma_2(c) = \tau_1^{A,B}(c)+i\tau_2^{A,B}(c) + \ln
  \epsilon^2,
  \]
  which proves the first statement of the proposition.
  
  Let us show now that for regular values of $c$ the 1-form
  $\tau_1(c)dc_1+\tau_2(c)dc_2$ is closed. For this we fix a regular
  value $c_0$ and introduce the following action integral, for $c$ in
  a small ball of regular values around $c_0$:
  \begin{equation}
    \label{equ:action}
    \mathcal{A}(c):=\int_{\gamma_c}\alpha,
  \end{equation}
  where $\alpha$ is any 1-form on some neighbourhood of $\Lambda_c$ in
  $M$ such that $d\alpha=\omega$ (which always exists since
  $\Lambda_c$ is Lagrangian), and $c\fleche\gamma_c$ is a smooth
  family of loops on the torus $\Lambda_c$ with the same homology
  class as the trajectory of the joint flow of $(H_1,H_2)$ during the
  time $(\tau_1(c),\tau_2(c))$. A simple argument (see for instance
  \cite[Lemma 3.6]{san-focus}) shows that
  $\deriv{\mathcal{A}(c)}{c_j}=\int_{\gamma_c}\kappa_j$, where
  $\kappa_j$ is the closed 1-form on $\Lambda_c$ defined by
  $\iota_{\ham{H_i}}\kappa=\delta_{i,j}$. In other words, the integral
  of $\kappa_j$ along a trajectory of the flow of $H_j$ measures the
  increase of the time $t_j$ along this trajectory. This means that
  \begin{equation}
    d\mathcal{A}(c) = \tau_1(c)dc_1+\tau_2(c)dc_2,
    \label{equ:daction}
  \end{equation}
  and thus proves the closedness of the right-hand side.
  
  Another way of proving this fact would be to apply the
  Liouville-Arnold theorem, which ensures that any 1-form
  $adc_1+bdc_2$, where $a$, $b$ depend smoothly on $c$ near a regular
  value, such that $(a,b)$ is in the \emph{period lattice}, is
  closed (see remark~\ref{rema:naturality}).
  
  Adding the fact that $\ln(c)dc$ is closed as a holomorphic 1-form,
  we obtain the closedness of $\sigma$ at any regular value of $c$,
  and hence at $c=0$ as well.
\end{demo}

\begin{rema}
  From this proposition, one easily recovers the result of
  \cite{zung2} stating that the monodromy of the Lagrangian fibration
  around a \ff\ fibre is generated by the matrix $\left(
    \begin{array}{cc}
      1 & 1 \\ 0 & 1
    \end{array}
  \right)$.
\end{rema}
Notice that the function $c\fleche\sigma_2(c)$ is defined modulo the
addition of a fixed constant in $2\pi\ZM$. We shall from now on assume
that $\sigma_2(0)\in[0,2\pi[$. This amounts to choosing the
determination of the complex logarithm in accordance with the
determination of $\tau_2$.
\begin{defi}
  \label{defi:invariant}
  Let $S$ be the unique smooth function defined in some neighbourhood
  of $0\in\RM^2$ such that $dS=\sigma$ and $S(0)=0$.  The Taylor
  expansion of $S$ at $c=0$ is called the \textbf{symplectic
    invariant} of Theorem \ref{theo:abstrait}. It is denoted by
  $(S)^\infty$.
\end{defi}
\begin{rema}
  \label{rema:action}
  Using equation~(\ref{equ:daction}), one can interpret $S$ as a
  \emph{regularised action integral}:
  \[
  S(c) = \mathcal{A}(c) - \mathcal{A}(0) + \re(c\ln c -c).
  \]
\end{rema}
\begin{rema}
  \label{rema:naturality}
  The formula \eqref{equ:daction} defines the 1-form
  $\tau=\tau_1dc_1+\tau_2dc_2$ independently of the choice of the
  coordinate system $(c_1,c_2)$.  Another (standard) way of viewing
  this is the following. First let $\B$ be the set of regular leaves
  of $\F$, and $\pi$ be the projection (which is a Lagrangian
  fibration) $\F\flechedroite{\pi}\B$.  The choice of a particular
  semi-global momentum map $F:=(H_1,H_2)$ for the system (near a
  Lagrangian leaf $\Lambda_c:=\pi^{-1}(c)$, for some $c\in\B$) is
  equivalent to the choice of a \emph{local chart} $\phi$ for $\B$
  near $c$: $F=\phi\circ\pi$.
  
  Then for each $c\in\B$, $T^*_c\B$ acts naturally on $\Lambda_c$ by
  the time-1 flows of the vector fields symplectically dual to the
  pull backs by $\pi$ of the 1-forms in $T^*_c\B$. This action extends
  to a Hamiltonian action in a neighbourhood of $\Lambda_c$ if and
  only if we restrict to closed 1-forms on $\B$. (In the local
  coordinates $(c_1,c_2)$ of $\B$ given by the choice of a momentum
  map $F=(H_1,H_2)$, the constant 1-forms $dc_1$, $dc_2$ act by the
  flows of $\ham{H_1}$, $\ham{H_2}$, respectively).
  
  The \emph{stabiliser} of this action form a particularly interesting
  lattice in $T^*_c\B$, which is another representation of the
  ``period lattice'' \cite{duistermaat}. It is the main point of the
  Liouville-Arnold theorem to show that, as $c$ varies, the points of
  this lattice are associated to \emph{closed} 1-forms, called
  \emph{period 1-forms}. (Indeed, in action-angle coordinates, the
  period 1-forms have constant coefficients). In our case, the period
  lattice is computed using a local chart given by Eliasson's theorem.
  First we see that this lattice has a privileged direction given by
  the $S^1$-action of $q_2$. Then we construct a ``minimal'' basis of
  this lattice by choosing the generator of this $S^1$-action (ie
  $2\pi dc_2$) together with the smallest transversal vector $\tau$
  that has positive coefficients on $dc_1$ and $dc_2$.  This is what
  we have done in this section.
\end{rema}

\section{Uniqueness}
\label{sec:uniqueness}
In order to show that the above invariant $(S)^\infty$ is indeed
symplectic and uniquely defined by the foliation, we need to prove
that it does not depend on any choice made to define them. A priori,
$(S)^\infty=(S)^\infty(\F,\chi)$ depends on the foliation $\F$ and on
the choice of the chart $\chi$ that puts a neighbourhood of the \ff\ 
point $m$ into normal form. It follows from the definition that if
$\phy$ is a symplectomorphism sending $\F$ to $\tilde{\F}$, then
$(S)^\infty(\tilde{\F},\tilde{\chi}) =
(S)^\infty(\F,\tilde{\chi}\circ\phy)$. So $(S)^\infty$ is well-defined
as a symplectic invariant of $\F$ if and only if, for any choice of
two chart $\chi$ and $\chi'$ putting a neighbourhood of $m$ into
normal form, $(S)^\infty(\F,\chi)=(S)^\infty(\F,\chi')$. This is
guaranteed by the following lemma:
\begin{lemm}
  \label{lemm:fonctionG}
  If $\phy$ is a local symplectomorphism of $(\RM^4,0)$ preserving the
  standard \ff\ foliation $\{q:=(q_1,q_2)=\textrm{const}\}$ near the
  origin, then there exists a unique germ of diffeomorphism
  $G:\RM^2\fleche\RM^2$ such that
  \begin{equation}
    \label{equ:fonctionG}
    q\circ \phy = G\circ q,
  \end{equation}
  and $G$ is of the form $G=(G_1,G_2)$, where
  $G_2(c_1,c_2)=\varepsilon_2 c_2$ and $G_1(c_1,c_2)-\varepsilon_1c_1$
  is flat at the origin, with $\varepsilon_j=\pm1$.
\end{lemm}
\begin{rema}
  This uniqueness statement about Eliasson's normal form does not
  appear in \cite{eliasson-these}.
\end{rema}
\begin{demo}[of the lemma]
  The existence of some unique $G$ satisfying \eqref{equ:fonctionG} is
  standard (because the leaves of the \ff\ foliation are locally
  connected around the origin). What interests us here are the last
  properties. As before, we use the complex coordinates
  $(z_1,z_2)\in\CM^2=\RM^4$, and $c=\bar{z_1}z_2\in\CM=\RM^2$. Let
  $\delta>0$ be such that $\phy$ is defined in the box
  $\mathcal{B}=\{|z_1|\leq 2\delta,|z_2|\leq 2\delta\}$.
   
  Since the flow of $q_2$ is $2\pi$-periodic, (\ref{equ:fonctionG})
  implies that the Hamiltonian vector field
  $\partial_1G_2\ham{q_1}+\partial_2G_2\ham{q_2}$ is also
  $2\pi$-periodic (with $2\pi$ as a primitive period). But on
  $\Lambda_0$ the only linear combinations of $\ham{q_1}$ and
  $\ham{q_2}$ that are periodic are the integer multiples of
  $\ham{q_2}$. Hence $\partial_1G_2(0)=0$ and $\partial_2G_2(0)=\pm
  1$.
  
  The flow of $q_1$ on $\Lambda_0$ is radial: any line segment $]0,A[$
  for some $A\in\Lambda_0$ is a trajectory. Then by
  (\ref{equ:fonctionG}) it image by $\phy$ must be a trajectory of
  $G_1\circ q$. Since $\phy$ is smooth at the origin, the image of
  $]0,A[$ for $A\in\mathcal{B}$ close enough to $0$ lies in some
  proper sector of the plane $\Pi\subset\Lambda_0$ containing
  $\phy(A)$ ($\Pi$ is either $\{z_1=0\}$ or $\{z_2=0\}$). But the only
  linear combinations of $\ham{q_1}$ and $\ham{q_2}$ which yield
  trajectories that are confined in a proper sector of $\Pi$ are the
  multiples of $\ham{q_1}$. Hence $\partial_2G_1(0)=0$. It follows now
  from the previous paragraph that $\partial_1G_1(0)\neq 0$ (since $G$
  is a local diffeomorphism).

  $\phy$ preserves the critical set of $q$; since left composition of
  $\phy$ by the symplectomorphism $(z_1,z_2)\fleche(-z_2,z_1)$ leaves
  \eqref{equ:fonctionG} unchanged (except for the sign of $G_1$), we
  may assume that each ``axis'' ($\{z_2=0\}$ and $\{z_1=0\}$
  respectively) is preserved by $\phy$. But then $\{z_2=0\}$ is the
  local unstable manifold for both $q_1$ and $G_1(q_1,q_2)$, which
  says that $\partial_1G_1(0)>0$.
  
  Using \eqref{equ:flowq1} and \eqref{equ:flowq2}, it is immediate to
  check that the joint flow of $(q_1,q_2)$ taken at the joint time
  $(-\ln|c/\delta|, \arg c)$ sends the point $(\bar{c},\delta)$ to the
  point $(\delta,c)$, and hence extends to a smooth and single-valued
  map $\Phi$ from a neighbourhood of $(0,\delta)$ to a neighbourhood
  of $(\delta,0)$.
  
  $\phy^{-1}\circ\Phi\circ\phy$ sends a neighbourhood of
  $\phy^{-1}(0,\delta)=(0,a)$ to a neighbourhood of
  $\phy^{-1}(\delta,0)=(b,0)$ and, because of
  \eqref{equ:fonctionG}, it is equal -- in the complement of the
  singular leaf $\Lambda_0$ -- to the joint flow of $G\circ q$ at the
  joint time $(-\ln|c/\delta|, \arg c)$, which is equal to the
  joint flow of $q$ at the joint time
  \[
  (-\partial_1G_1\ln|c/\delta|+\partial_1G_2\arg c,
  -\partial_2G_1\ln|c/\delta|+\partial_2G_2\arg c).
  \] 
  Since $\phy^{-1}\circ\Phi\circ\phy$ is smooth at the origin, we
  obtain by restricting the first component of this map to the
  ``Poincaré'' surface $\{(\bar{c},a)$ with $c$ near $0$ in $\CM\}$
  that the map:
  \begin{equation}
    c\fleche \exp\left((1-\partial_1G_1)\ln|c|+\partial_1G_2\arg
      c+i\left((\partial_2G_2-1)\arg c
        -\partial_2G_1\ln|c|\right)\right)
    \label{equ:flow}
  \end{equation}
  is single-valued and smooth at the origin. (We have factored out the
  terms $\exp(\partial_jG_1\ln\delta)$, $j=1,2$, which are obviously
  smooth.)

  The single-valuedness of \eqref{equ:flow} implies that
  $\partial_1G_2\equiv 0$ and $\partial_2G_2\in\ZM$. Hence
  $\partial_2G_2=\pm 1$.
  
  Now the smoothness of \eqref{equ:flow} says that the following two
  functions:
  \[
  c\fleche (1-\partial_1G_1)\ln|c| \qquad \textrm{ and } \qquad
  c\fleche -\partial_2G_1\ln|c|
  \]
  are smooth at the origin, which easily implies that
  $(1-\partial_1G_1)$ and $\partial_2G_1$ are flat at the origin,
  yielding the result.
\end{demo}
Suppose we define two semi-global invariants $(S)^\infty(\F,\chi)$ and
$(S)^\infty(\F,\tilde{\chi})$ by choosing two different charts $\chi$
and $\tilde{\chi}$ which put a neighbourhood of the \ff\ point into
normal form. As before, one defines the momentum maps $F$ and
$\tilde{F}$, which are the extensions to $\F$ of $q\circ\chi$ and
$q\circ\tilde{\chi}$, and computes the corresponding period 1-forms
$\tau$ and $\tilde{\tau}$.  Then we can invoke the lemma to
$\phy=\tilde{\chi}\circ \chi^{-1}$, and the conclusions apply to
$G=\tilde{F}F^{-1}$.
 
Suppose that $\varepsilon_j=1$, $j=1,2$, i.e. $G$ is infinitely
tangent to the identity.  Then the same type of arguments as above (a
logarithm cannot compete against a flat term) shows that, since the
vector fields $\ham{H_j}$ and $\ham{\tilde{H}_j}$ are infinitely
tangent to each other, $\tau$ and $\tilde{\tau}$ must differ by a flat
term. Actually, since by remark \ref{rema:naturality}
$G^*\tilde{\tau}$ is also a period 1-form associated with the momentum
map $F$, one has $\tau=G^*\tilde{\tau}$.  This implies that
$\sigma(c)=\tau(c)+\Re(\ln c dc)$ and
$\tilde{\sigma}=(G^{-1})^*\sigma$ differ by a flat form, hence
$(S)^\infty(\F,\chi)= (S)^\infty(\F,\tilde{\chi})$.

If $\varepsilon_2=-1$, it suffices to compose with the
symplectomorphism $(x,\xi)\fleche (-x,-\xi)$, which sends $(q_1,q_2)$
to $(q_1,-q_2)$ and leaves $\sigma$ invariant (both $\sigma_2$ and
$dc_2$ change sign). An analogous remark holds with the
symplectomorphism $(z_1,z_2)\fleche(-z_2,z_1)$, which sends
$(q_1,q_2)$ to $(-q_1,q_2)$ and leaves $\sigma$ invariant, while
changing the sign of $\varepsilon_1$.

\section{Injectivity}

Let $\F$ and $\tilde{\F}$ are two singular foliations of simple \ff\ 
type on the symplectic manifolds $(M,\omega)$ and
$(\tilde{M},\tilde{\omega})$. Assume that they have the same invariant
$(S)^\infty(\F)=(S)^\infty(\tilde{\F})\in\RM\formel{X,Y}_0$. We shall
prove here that $\F$ and $\tilde{\F}$ are semi-globally equivalent,
ie.  there exists a foliation preserving symplectomorphism between
some neighbourhoods of the \ff\ leaves.

For each of the foliations $\F$ and $\tilde{\F}$, we choose a chart of
Eliasson's type around the \ff\ point, and thus define the period
1-forms $\tau$ and $\tilde{\tau}$ on $(\RM^2\setminus\{0\},0)$. The
hypothesis implies that there is a smooth closed 1-form
$\pi=\pi_1dc_1+\pi_2dc_2$ on $(\RM^2,0)$ whose coefficients are flat
functions of $c$ at the origin such that
\[
\tilde{\tau} = \tau + \pi.
\]
\begin{lemm}
  One can chose symplectic charts of Eliasson's type at the \ff\ 
  points in such a way that $\pi=0$, ie:
  \[
  \tilde{\tau} = \tau.
  \]
\end{lemm}
\begin{demo}
  1. We first prove that there exists a local diffeomorphism $G$ of
  $(\RM^2,0)$ isotopic to the identity such that $(G^{-1})^*\tau =
  \tilde{\tau}$. We wish to realise $G$ as $G_1$ where $G_t$ is a flow
  satisfying
  \[
  G_t^*(\tau+t\pi)=\tau.
  \]
  This amounts to finding the associated vector field $Y_t$ which must
  satisfy \[ d(\iota_{Y_t}(\tau+t\pi)) = -\pi.
  \]
  We can write $\pi=dP$ for some smooth function $P$ which is flat at
  $0$. Assume we look for a field $Y_t$ of the form
  $Y_t=f_t(c)\deriv{}{c_1}$. We obtain the following equation:
  \[
  f_t(c)=\frac{-P(c)}{\tau_1(c)+t\pi_1} = \frac{-P(c)}{\ln|c|
    -\sigma_1(c)+t\pi_1}.
  \]
  Since $P$ is flat at $0$, the right-hand-side is indeed a (flat)
  smooth function depending smoothly on $t$, and the result is proved.

  2. Notice also that $G$ is infinitely tangent to the identity, and
  moreover leaves the second variable $c_2$ unchanged. Now we show
  that for any diffeomorphism $G$ of $(\RM^2,0)$ sharing these
  properties (which are those of Lemma \ref{lemm:fonctionG}) there
  exists a symplectomorphism $\chi$ near the \ff\ point $m$ such that
  \[
  G(q_1,q_2)\circ\chi = (q_1,q_2).
  \]
  Here again we seek $\chi$ as the time-1 map of the flow of some
  vector field $X_t$. Of course we shall look now for a Hamiltonian
  vector field $X_t=\ham{f_t}$ to ensure the symplecticity of
  $\chi_t$. Then the requirement
  \[
  \chi_t^* q_t = q_0,
  \]
  where $q_t=(q_{t,1},q_{t,2})\egdef tG(q_1,q_2)+(1-t)(q_1,q_2)$,
  leads to the following system
\[
\begin{array}{lll}
\{f_t,q_{t,1}\} & = & g_1 \\
\{f_t,q_{t,2}\} & = & 0,
\end{array}
\]
with $(g_1,0)=(q_1,q_2)-G(q_1,q_2)$. By hypothesis $g_1$ is a flat
function at the origin, and the fact that $\{q_{t,1},q_{t,2}\}\equiv
0$ implies that $\{g_1,q_{t,2}\}=0$.  Moreover the quadratic part of
$q_t$ is $q_0$, so we know (see \cite{eliasson-these}) that such a
system admits a solution $f_t$.

It remains to put all our remarks together: Point 2) shows that left
composition by $\chi$ of the Eliasson chart we have chosen at $m$ is
again an admissible chart of Eliasson's type, yielding the new
momentum map $G(q_1,q_2)$. Using the $G$ obtained at Point 1) and in
view of the naturality property (remark \ref{rema:naturality}), the
new period 1-form (denoted by $\tau$ again) satisfies
$\tau=\tilde{\tau}$.
\end{demo}

We are now is position to construct the required equivalence.
Applying the lemma we get a local symplectomorphism that allows us to
identify some neighbourhoods $U$ and $\tilde{U}$ of the \ff\ points
$m$ and $\tilde{m}$, and two momentum maps $F$ and $\tilde{F}$ (both
equal to $(q_1,q_2)$ inside their respective neighbourhoods of the
\ff\ points) which define the same closed 1-form $\sigma$ on
$(\RM^2,0)$. We denote $\Lambda_c=F^{-1}(c)$ and
$\tilde{\Lambda}_c=\tilde{F}^{-1}(c)$.

Let $\mathfrak{U}$ be an open ball strictly contained in $U$, let
$\Sigma_u\subset\mathfrak{U}$ be a transversal section as defined in
the proof of Proposition \ref{prop:period}, and construct in the same
way $\tilde{\Sigma}_u$ for the foliation $\tilde{\F}$ (so that
$\Sigma_u$ and $\tilde{\Sigma}_u$ are identified by the above
symplectomorphism).  Reduce $\F$ (and $\tilde{\F}$) to the
neighbourhoods of the \ff\ leaves composed of the leaves intersecting
$\Sigma_u$ (or $\tilde{\Sigma}_u$). We construct our equivalence by
extending the identity outside $\mathfrak{U}$. Let
$x\in\Lambda_c\setminus\mathfrak{U}$, and define
$t(x)\in]0,\tau_1(c)[$ to be the smallest time it takes for the point
$\Sigma_u\cap\Lambda_c$ to reach the $\ham{H_2}$-orbit of $x$. (Recall
that $H_2$ generates an $S^1$ action.) Now define $s(x)\in\RM/2\pi\ZM$
as the remaining time to finally reach $x$ under the $\ham{H_2}$-flow.
To this $x$ we associate the point $\tilde{x}\in\tilde{\F}$ obtained
from the point $\tilde{\Sigma}_u\cap\tilde{\Lambda}_c$ by letting the
joint flow of $\tilde{F}$ act during the times $(t(x),s(x))$. This map
--- let's call it $\Psi$ --- is well defined because of the equality
$\tau=\tilde{\tau}$.  It is a bijection since the inverse is equally
well-defined just by interchanging the roles of $\F$ and $\tilde{\F}$.
Between $U$ and $\tilde{U}$, $\Psi$ is a symplectomorphism since
through Eliasson's charts, it is just the identity. Concerning now the
symplecticity of $\Psi$ in the complement of the singular points, one
can prove it for $c\neq 0$ (which is sufficient by continuity) by
invoking the Liouville-Arnold theorem, which shows that $\Psi$ is
symplectically conjugate to a translation in the fibres. Then the
symplectic property near the singular points implies that this
translation must be symplectic everywhere. A similar argument using
the less sophisticated Darboux-Carathéodory theorem could also do. But
the simplest is maybe the following. It is clear from the construction
that $\Psi$ is equivariant with respect to the joint flows of our
Hamiltonian dynamics:
\begin{equation}
  \forall (t_1,t_2), \qquad \Psi\circ \phy_{t_1,t_2}  =
  \tilde{\phy}_{t_1,t_2}\circ \Psi,
  \label{equ:equivariant}
\end{equation}
where $\phy_{t_1,t_2}$ and $\tilde{\phy}_{t_1,t_2}$ are the joint
flows of $F$ and $\tilde{F}$ at the joint time $(t_1,t_2)$.
Using~\eqref{equ:equivariant} together with the fact that
$\tilde{\phy}_{t_1,t_2}$ is symplectic, we see that
$\phy_{t_1,t_2}^*(\Psi^*\tilde{\omega})=\Psi^*\tilde{\omega}$; in
other words, $\Psi^*\tilde{\omega}$ is invariant under the joint flow
$\phy_{t_1,t_2}$. Since $\omega$ has the same property, so has
$\Psi^*\tilde{\omega}-\omega$. Since $\Psi^*\tilde{\omega}-\omega=0$
near $m$, it must vanish as well on the whole $\F$.

\section{Surjectivity}
\label{sec:surjectivity}
We prove here that any formal power series $(S)^\infty\in
\RM\formel{X,Y}_0$ is the symplectic invariant --- in the sense of
Definition \ref{defi:invariant} --- of some Liouville foliation of
simple \ff\ type. More precisely, we construct a foliation $\F$
together with a local chart $\chi$ that puts a neighbourhood of the
\ff\ point into normal form such that (using the notation of
Section~\ref{sec:uniqueness}) $(S)^\infty=(S)^\infty(\F,\chi)$.
Another proof of this result has been proposed by Castano-Bernard
\cite{castano-proof}.

Using the same notations as before, we let $(q_1,q_2)=\bar{z_1}z_2$ be
the standard \ff\ fibration $\RM^4\simeq\CM^2\fleche\CM\simeq\RM^2$
defined in (\ref{equ:cartan}). The joint flow will be denoted by
$\phy_{t_1,t_2}$.

Invoking Borel's construction, let $S\in\Cinf(\RM^2)$ be a function
vanishing at the origin and whose Taylor series is $(S)^\infty$.  We
shall denote by $S_1$, $S_2$ the partial derivatives $\partial_XS$ and
$\partial_YS$, respectively.

Let us define two ``Poincaré'' surfaces in $\CM^2$ by means of the
following embeddings of the ball
$D_\varepsilon=B(0,\varepsilon)\subset\CM$, for some
$\varepsilon\in]0,1[$:
\[
\begin{array}{l}
\Pi_1(c)=(\bar{c},1)\\
\Pi_2(c)=(e^{S_1(c)+iS_2(c)},ce^{-S_1(c)+iS_2(c)}).
\end{array}
\]
Notice that for each $c$, the points $\Pi_j(c)$, $j=1,2$ belong to the
(non-compact) Lagrangian submanifold $\Lambda_c:=\{\bar{z_1}z_2=c\}$.
$\Pi_j(D_\varepsilon)$, $j=1,2$ are smooth 2-dimensional manifolds
constructed in such a way that for any $c\neq 0$, $\Pi_2(c)$ is the
image of $\Pi_1(c)$ by the joint flow of $(q_1,q_2)$ at the time
$(S_1(c)-\ln|c|,S_2(c)+\arg(c))$.

Let $\Phi$ be this diffeomorphism, defined on all
$\Pi_j(D_\varepsilon)$ by the embeddings:
\[
\vcenter{\xymatrix{
    \Pi_1(D_\varepsilon) \ar[rr]^{\Phi} & & \Pi_2(D_\varepsilon) \\
    & D_\varepsilon\ar[ul]^{\Pi_1}\ar[ur]_{\Pi_2} & }}
\]

$\Pi_1(D_\varepsilon)$ and $\Pi_2(D_\varepsilon)$ are transversal to
the Lagrangian foliation, and $\Phi$ can be extended uniquely to a
diffeomorphism between small neighbourhoods of $\Pi_1(D_\varepsilon)$
and $\Pi_2(D_\varepsilon)$ by requiring that it commute with the joint
flow:
\begin{equation}
  \Phi\left(\phy_{t_1,t_2}(m)\right) =
  \phy_{t_1,t_2}\left(\Phi(m)\right).
  \label{equ:commute}
\end{equation}

\begin{lemm}
  $\Phi$ is a symplectomorphism.
\end{lemm}
\begin{demo}
  One can write $\Phi$ in terms of $\Pi_1$ and $\Pi_2$ and check the
  result by explicit calculation.  However, the reason why it works is
  the following:
  
  Since we already know that $\Phi$ is smooth, it is enough to prove
  the lemma outside of the singular Lagrangian $\Lambda_0$.  So fix
  $c_0\neq 0$; we can construct a Darboux-Carathéodory chart
  $(x,\xi)\in\RM^4$ in a connected open subset of $\Lambda_{c_0}$
  containing both $\Pi_1(c_0)$ and $\Pi_2(c_0)$.  In these
  coordinates, the momentum map is $(\xi_1,\xi_2)$ and the flow is
  linear: $\phy_{t_1,t_2}$ is the translation by $(t_1,t_2)$ in the
  $x$ variables.
  
  Through this chart, $\Phi$ is by construction a ``fibre
  translation'':
  \begin{equation}
    \Phi(x,\xi) = (x+f(\xi),\xi),
    \label{equ:translation}
  \end{equation}
  where
  \begin{equation}
    \label{equ:f}
    f(\xi)=(S_1(\xi),S_2(\xi))+ (\ln|\xi|,-\arg(\xi)).
  \end{equation}
  
  Now, it is easy to check that (\ref{equ:translation}) defines a
  symplectomorphism if and only if the 1-form
  \[ 
  f_1(\xi)d\xi_1+f_2(\xi)d\xi_2
  \]
  is closed.  In our case the closedness is automatic since
  $S_1dX+S_2dY=dS$.
\end{demo}
Let $\Sigma_j$, $j=1,2$ be the $S^1$-orbit of $\Pi_j(D_\varepsilon)$.
Construct a 4-dimensional cylinder $\mathcal{C}$ by letting the
$q_2$-flow take $\Sigma_1$ to $\Sigma_2$, namely:
\[
\mathcal{C} := \overline{\bigsqcup_{c\in D_\varepsilon\setminus\{0\}}
  \mathcal{C}_c}
\]
where $\mathcal{C}_c\subset\Lambda_c$ is the 2-dimensional cylinder
spanned by $\phy_{t_1,t_2}(\Pi_1(c))$, for $(t_1,t_2)\in
[0,S_1(c)+\ln|c|]\times[0,2\pi]$.
\begin{figure}[htbp]
  \begin{center}
    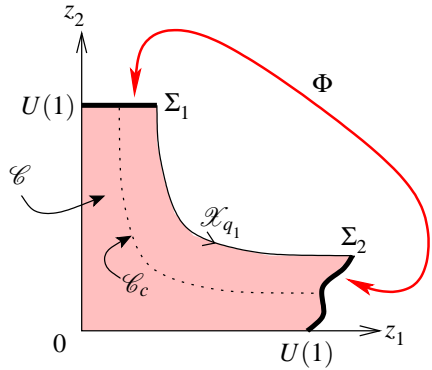
    \caption{Construction of the symplectic manifold $M$}
    \label{fig:construction}
  \end{center}
\end{figure}
Finally, let $M$ be the symplectic manifold obtained by gluing the two
ends $\Sigma_j$ of the cylinder $\mathcal{C}$ using the
symplectomorphism $\Phi$. Since $\Phi$ preserves the momentum map
$(q_1,q_2)$, the latter yields a valid momentum map $F$ on $M$. The
corresponding Lagrangian foliation $F^{-1}(c)$ is given by
$\mathcal{C}_c$ with its two ends identified by $\Phi$. In particular
all leaves are compact and the foliation is of simple \ff\ type.

The $S^1$ action is unchanged, while the transversal period
$(\tau_1(c),\tau_2(c))$ on $F^{-1}(c)$ is by construction the time it
takes for the joint flow to reach $\Pi_2(c)$ from $\Pi_1(c)$, ie
\[
(\tau_1(c),\tau_2(c)) = (S_1(c)-\ln|c|,S_2(c)+\arg(c)).
\]
Then by definition \ref{defi:invariant} the symplectic invariant of
the foliation is given by the Taylor expansion of the primitive of the
1-form $S_1dc_1+S_2dc_2$ vanishing at $0$, ie. $(S)^{\infty}$.

\section{Further remarks}

\paragraph{Multiple \ff.} Assume now that the singular fibre
$\Lambda_0$ carries $k$ \ff\ points $m_0, \dots, m_{k-1}$. Then
$\Lambda_0$ is a $k$-times pinched torus, and
Theorem~\ref{theo:abstrait} can be generalised. In this case, the
regularisation of the action integral $S$ must take into account all
the singular points. In order to do this, one has to consider $k-1$
local invariants, which are also formal power series in
$\RM\formel{X,Y}$, and which measure the obstruction to construct a
semi-global momentum map that is in Eliasson normal form
simultaneously at two different singular points. Here follows a sketch
of the argument.

Let $F$ be a semi-global momentum map. At each point $m_j$ one has a
local normal form $F\circ\phy_j=G^j(q_1,q_2)$. Because of Lemma
\ref{lemm:fonctionG}, one can extend $q_2$ to a periodic Hamiltonian
on a whole neighbourhood of $\Lambda_0$, and one can always assume
that $\phy_j$ is orientation preserving --- that means we fix once and
for all the sign of the $\varepsilon_j$. If now $F$ if of the form
$(H_1,q_2)$ then $G^j$ takes the form
$G^j(q_1,q_2)=(F^j(q_1,q_2),q_2)$.  By the implicit function theorem,
$F^j$ is locally invertible with respect to the variable $q_1$. Let
$(F^j)^{-1}$ be this inverse, and define $G^{i,j}=(F^i)^{-1}F^j$.
Again by Lemma \ref{lemm:fonctionG}, the Taylor expansions of
$G^{i,j}$ are invariants of the foliation.

Assume the points $m_i$ are ordered according to the flow of $H_1$,
with indices $i\in\ZM/k\ZM$.
Similarly to the case $k=1$, one can define a regularised period
1-form $\sigma$ by the following formula:
\begin{equation}
  \sigma:=\sum_{i=0}^{k-1}(G_0^{-1}G_i)^*\left(\sigma_1^{i,i+1}(c)dc_1 +
  \sigma_2^{i,i+1}(c)dc_2\right), 
  \label{equ:periode_multiple}
\end{equation}
with
\begin{equation}
  \left\{
    \begin{array}{ccl}
      \sigma_1^{i,i+1}(c) & = & \tau_1^{i,i+1}(c)+\Re(\ln c) \\
      \sigma_2^{i,i+1}(c) & = & \tau_2^{i,i+1}(c)-\Im(\ln c)
    \end{array}
  \right.,
  \label{equ:cocycle}
\end{equation}
where $(\tau_1^{i,i+1}(c), \tau_2^{i,i+1}(c))$ are the smallest
positive times needed to reach $A_{i+1}(c)$ from $A_i(c)$ under the
flow of $(G^i)^{-1}\circ F$ --- which is the momentum map $(q_1,q_2)$
in the normal form coordinates near point $A_i$. Here we have chosen a
point $A_i(c)$ in a Poincaré section of each local stable manifold
near $m_i$. Of course $\sigma_j^{i,i+1}(c)$ depends heavily on the
choice of $A_i$ and $A_{i+1}$, but the sums appearing in
\eqref{equ:periode_multiple} does not, and the resulting 1-form
$\sigma$ is closed.
\begin{figure}[htbp]
  \begin{center}
    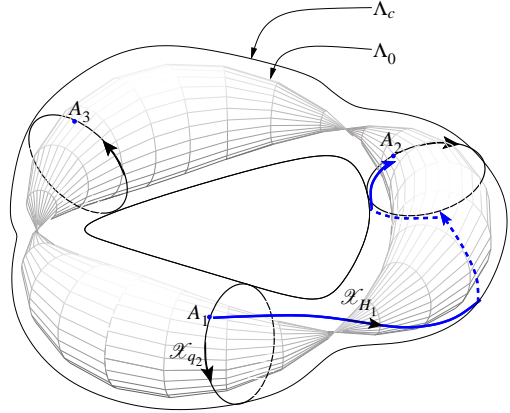
    \caption{The multi-pinched torus}
    \label{fig:multiple}
  \end{center}
\end{figure}
Notice that the definition of $\sigma$ depends on the choice of a
start point $m_0$. Thus we are here classifying a singular foliation
with a distinguished \ff\ point $m_0$.

\begin{figure}[htbp]
  \begin{center}
    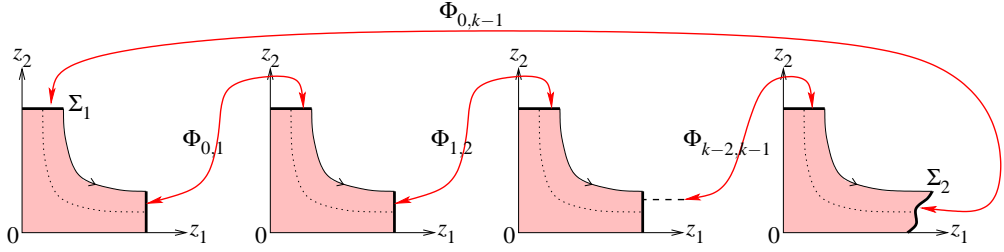
    \caption{multiple gluing}
    \label{fig:constructionX}
  \end{center}
\end{figure}

Let $(S)^\infty$ be the Taylor series of the primitive of $\sigma$
vanishing at the origin.  Then $(S)^\infty$ and the $k-1$ ordered
invariants $(G^{i,i+1})^\infty$ are independent and entirely classify
a neighbourhood of the critical fibre $\Lambda_0$ with distinguished
point $m_0$. The arguments of the proof are similar to the ones of the
case $k=1$. An abstract construction of a foliation admitting a given
set of invariants is proposed in Figure~\ref{fig:constructionX}. There
the local pictures are described by canonical coordinates respectively
given by $(q_1,q_2)$, $(G^{1,2}(q_1,q_2),q_2)$,
$(G^{1,2}(G^{2,3}(q_1,q_2),q_2),q_2)$, etc. and the gluing
diffeomorphisms $\Phi_{i,i+1}$ are constructed as in
section~\ref{sec:surjectivity} using the following functions,
respectively: $S_{0,1}=S_{1,2}=\cdots=S_{k-2,k-1}=0$ and $S_{0,k-1}$
is a resummation of $(S)^\infty$.
\begin{rema}
  We can regard the reduced space $\Lambda_0/S^1$ as a cyclic graph
  $\mathcal{G}$ whose vertices are the \ff\ points $m_i$, and which is
  oriented by the flow of $H_1$. For each edge $[i,i+1]$ one can
  define a 1-form
  \[
  \sigma^{i,i+1}:=(G_0^{-1}G_i)^*\left(\sigma^{i,i+1}_1dc_1+
    \sigma^{i,i+1}_2dc_2\right) \in \Omega^1(D)
  \] 
  (for some fixed small disc $D$ around the origin in $\RM^2$). This
  defines a 1-cocycle on $\mathcal{G}$ with values in the vector space
  $\Omega^1(D)$. If one varies the points $A_j$, this cocycle is
  easily seen to change by a coboundary; hence the set of
  $\{\sigma^{i,i+1}\}$ naturally defines a well-defined cohomology
  class on $\mathcal{G}$. Be the same argument as in the case $k=1$
  (ie. essentially Arnold-Liouville's theorem) this class is
  \emph{closed}, in the sense that the cochain $\{\sigma^{i,i+1}\}$,
  modulo some coboundary, can be chosen to consist only of closed
  1-forms. Hence we end up with a class $[\sigma]\in
  H^1(\mathcal{G},H^1(D))$. Since $\mathcal{G}$ is homeomorphic to a
  circle, $H^1(\mathcal{G},H^1(D))\simeq H^1(D)$ and $[\sigma]$ is
  represented by the de Rham cohomology class of the closed 1-form
  $\sigma=\sum\sigma^{i,i+1}$ defined in~\eqref{equ:periode_multiple}.
  
  Now, the functor that produces Taylor series of 1-forms can be
  applied to the coefficients of this cochain, yielding a cocycle with
  values in formal closed 1-forms and whose class is represented by
  the differential of our invariant $(S)^\infty$.
\end{rema}

\paragraph{``Exact'' version.}
If one intends to extend the results to a semiclassical setting,
general symplectomorphism do not suffice: one needs to control the
action integrals (in the standard semiclassical pseudo-differential
theory, a potential $\alpha$ for the symplectic form: $d\alpha=\omega$
is part of the data). In view of Remark~\ref{rema:action}, this is
naturally done by including the constant term in the Taylor series of
$S$ as being the integral
\[
S_0:=\int_{\gamma_0}\alpha,
\]
where $\gamma_0$ is the generator of $H_1(\Lambda_0)$.

\paragraph{Acknowledgements.} This article answers a question that
J.J. Duistermaat asked me when I was in Utrecht in 1998. I wrote then
a short --- and incomplete --- draft, and that was it. Two years after
I finally read the note by Toulet and al. \cite{dufour-mol-toul}, and
a fruitful discussion with Richard Cushman made me realise that I had
the result at hand. I wish to thank him for this. I would also like to
thank Ricardo Castano-Bernard for interesting discussions, and for
showing me an alternative proof of the ``surjectivity'' part using
general arguments developed for mirror symmetry via special Lagrangian
fibrations.

After I wrote this article, P. Molino informed me of an unpublished
work of his (in collaboration with one of his students
\cite{grossi-molino}) concerning the same problem. They defined a
similar invariant, and the classification result was conjectured.

\bibliographystyle{plain}
  
\bibliography{bibli}

\begin{thebibliography}{10}

\bibitem{bolsinov-fomenko}
A.V. Bolsinov and A.T. Fomenko.
\newblock Orbital equivalence of integrable hamiltonian systems with two
  degrees of freedom. a classification theorem {I}.
\newblock {\em Russian Acad. Sci. Sb. Math.}, 81(2):421--465, 1995.
\newblock (translated from the Russian).

\bibitem{castano-proof}
R.~Castano-Bernard, 2000.
\newblock personal communication.

\bibitem{dufour-mol-toul}
J.-P. Dufour, P.~Molino, and A.~Toulet.
\newblock Classification des syst{\`e}mes int{\'e}grables en dimension 2 et
  invariants des mod{\`e}les de {F}omenko.
\newblock {\em C. R. Acad. Sci. Paris S{\'e}r. I Math.}, 318:949--952, 1994.

\bibitem{duistermaat}
J.J. Duistermaat.
\newblock On global action-angle variables.
\newblock {\em Comm. Pure Appl. Math.}, 33:687--706, 1980.

\bibitem{eliasson-these}
L.H. Eliasson.
\newblock {\em Hamiltonian systems with {P}oisson commuting integrals}.
\newblock PhD thesis, University of Stockholm, 1984.

\bibitem{eliasson}
L.H. Eliasson.
\newblock Normal forms for hamiltonian systems with {P}oisson commuting
  integrals -- elliptic case.
\newblock {\em Comment. Math. Helv.}, 65:4--35, 1990.

\bibitem{grossi-molino}
D.~Grossi.
\newblock Syst{\`e}mes int{\'e}grables non d{\'e}g{\'e}n{\'e}r{\'e}s ---
  invariants symplectiques d'une singularit{\'e} focus-focus simple.
\newblock M{\'e}moire de {DEA}, universit{\'e} de Montpellier, 1995.

\bibitem{zung}
Z.~Nguy{\^e}n~Ti{\^e}n.
\newblock A topological classification of integrable hamiltonian systems.
\newblock In R.~Brouzet, editor, {\em S{\'e}minaire Gaston Darboux de
  g{\'e}ometrie et topologie diff{\'e}rentielle}, pages 43--54. Universit{\'e}
  Montpellier II, 1994-1995.

\bibitem{zung2}
Z.~Nguy{\^e}n~Ti{\^e}n.
\newblock A note on focus-focus singularities.
\newblock {\em Diff. Geom. Appl.}, 7(2):123--130, 1997.

\bibitem{san-focus}
S.~V{\~u}~Ng{\d o}c.
\newblock Bohr-{S}ommerfeld conditions for integrable systems with critical
  manifolds of focus-focus type.
\newblock {\em Comm. Pure Appl. Math.}, 53(2):143--217, 2000.

\bibitem{san-fn}
S.~V{\~u}~Ng{\d o}c.
\newblock Formes normales semi-classiques des syst{\`e}mes compl{\`e}tement
  int{\'e}grables au voisinage d'un point critique de l'application moment.
\newblock {\em Asymptotic Analysis}, 24(3,4):319--342, 2000.

\bibitem{weinstein-symplectic}
A.~Weinstein.
\newblock Symplectic manifolds and their lagrangian submanifolds.
\newblock {\em Adv. in Math.}, 6:329--346, 1971.

\end{thebibliography}

\end{document}